# A modeler's guide to handle complexity in energy systems optimization

Leander Kotzur*,[a,f] Lars Nolting,[d,f] Maximilian Hoffmann,[a,f] Theresa Groß,[a,f] Andreas Smolenko,[c,f] Jan Priesmann,[d,f] Henrik Büsing,[c,f] Robin Beer,[a,f] Felix Kullmann,[a,f] Bismark Singh,[d] Aaron Praktiknjo,[d,f] Detlef Stolten[a,e,f] and Martin Robinius[a,f]

**Abstract**

Determining environmentally- and economically-optimal energy systems designs and operations is complex. In particular, the integration of weather-dependent renewable energy technologies into energy system optimization models presents new challenges to computational tractability that cannot only be solved by advancements in computational resources. In consequence, energy system modelers must tackle the complexity of their models by applying various methods to manipulate the underlying data and model structure, with the ultimate goal of finding optimal solutions. As which complexity reduction method is suitable for which research question is often unclear, herein we review different approaches for handling complexity. We first analyze the determinants of complexity and note that many drivers of complexity could be avoided a priori with a tailored model design. Second, we conduct a review of systematic complexity reduction methods for energy system optimization models, which can range from simple linearization performed by modelers to sophisticated multi-level approaches combining aggregation and decomposition methods. Based on this overview, we develop a guide for energy system modelers who encounter computational limitations.

Keywords: Energy system optimization, MILP, LP, decomposition, capacity expansion, aggregation

## 1 Introduction

### 1.1 Energy systems optimization

The design and operation of energy systems with minimal environmental and economic impacts is highly complex, as energy supply and demand must be spatially- and temporally-balanced, with an ever-increasing set of generation units, storage technologies, transmission options, and load management alternatives. The analytical solving of these problems is no longer feasible, and instead requires the use of mathematical energy system optimization models (ESOMs) to identify the optimal design and operation [1].

The theoretical limitation to the application of such optimization models is our ability to structure them in the form of a mathematical program [2]. These programs have a broad range of applications; for instance, to determine train routes and schedules [3], production planning [4], emergency logistics [5] or, as previously mentioned, energy systems [6], while limited resources and time stress efficient solution processes of those [7].

For the case of the energy sector, the concept to determine necessary future capacities using simple scenario-based models dates back to the 1950s [8]. Simultaneously, the first concepts for achieving profitability of energy trades such as peak-load pricing [9-11] have their origin in the same decade. During the 1960s and 1970s, the rapidly growing energy demand, as well as an advancing liberation of the energy market [5, 6], drove the development of more complex models in order to maintain the security of supply at every point in time and stay profitable despite an increasing number of market competitors. As a result, the first optimization-based models such as BESOM [12] emerged during this period. During the 1980s, environmental awareness and tech-

[a.] Institute of Techno-economic Systems Analysis (IEK-3), Forschungszentrum Jülich GmbH, Wilhelm-Johnen-Str., D-52428, Germany
[b.] Friedrich-Alexander-Universität Erlangen Nürnberg (FAU), Department of Mathematics, Cauerstr. 11, 91058 Erlangen, Germany
[c.] Institute for Advanced Simulation – Jülich Supercomputing Centre (IAS-JSC), Forschungszentrum Jülich GmbH, Wilhelm-Johnen-Str., D-52428, Germany
[d.] Chair for Energy System Economics (FCN-ESE), Institute for Future Energy Consumer Needs and Behavior, E.ON Energy Research Center, RWTH Aachen University, Mathieustrasse 10, 52074 Aachen, Germany
[e.] Chair of fuel cells, RWTH Aachen University, c/o Institute of Techno-economic Systems Analysis (IEK-3), Forschungszentrum Jülich GmbH, Wilhelm-Johnen-Str., D-52428, Germany
[f.] JARA-ENERGY, Jülich Aachen Research Alliance, 52425 Jülich, Germany
* corresponding author: l.kotzur@fz-juelich.de

nological advances first led to the consideration of renewable energy sources [13]. As many renewable energy sources are both, intermittent and non-dispatchable, their consideration led to a third major development of energy system modeling, namely the consideration of a finite number of regions as well as multiple discrete time steps in order to capture various demand and supply situations.

In the realm of energy system their application can be distinguished between top-down ESOMs approaches that address economic, political, and social aspects, which directly drive the evolution of energy systems [14], and bottom-up ESOMs that focus on detailed technological modelling and specific system design. Well-known representatives of this group are e.g. LEAP [15], EFOM [16], BESOM [12], MARKAL [17], MESSAGE [14, 18], IKARUS [6], PERSEUS [19], TIMES [20-23] and, recently, DESOD [24], DER-CAM [25], CALLIOPE [26], OEMOF [27], URBS [28], PYPSA [29] and FINE [30], while an overview of their type and application can be found in Appendix Table 3. A broader review of the scope and applications of bottom-up ESOMs is given by Groissböck [31], Lopion et al. [6], and Ringkjøb et al. [32].

**1.2 Increasing energy system complexity**

The quality and availability of the data required to parameterize these models are steadily improving, but the amount of input data required directly impacts the size of the related optimization problem and with it the requirement for processing resources and finding an optimal solution within a reasonable timeframe [33]. The non-linearities of the objective or constraint functions, or a large number of system variables and uncertainties, can even risk the identification of a feasible solution [34].

In particular, this makes the integration of renewable energy technologies into the bottom-up ESOMs challenging due to their often pervasive nonlinear structures and an increased requirement for spatiotemporal resolution [35-38]. In addition to the required increase in the granularity of the ESOMs, an increase in the connectivity of the different energy sectors must be considered, which is often referred to as 'sector-coupling' [39-41]. As an example, the electricity sector becomes strongly linked to the heat sector via heat pumps and other power-to-heat technologies. This coupling enables the construction of low-carbon energy systems with renewable electricity supply, but increases the intricacy of these systems and significantly threatens their viability, let alone the solvability of their corresponding models [32]. Thereby, data availability is improving, allowing for much more accurate model representations. E.g., hourly weather data is available for multiple decades across continental scales which represents huge amounts of data and decision which cannot be properly managed by simple mathematical programs.

Thus, settling for sub-optimal or merely feasible energy system design solutions will be the rule rather than the exception if there is no development towards better computational tractability to solve energy system models within an acceptable timeframe.

**1.3 Development of computational resources**

While one would expect that no major improvements in the solving of algorithms would be required in light of Moore's Law and the corresponding increasing availability of computer resources over the last few decades [42], the computational tractability of mathematical programs remains greatly limited [43].

This is caused by resource development: According to Moore's Law, the number of transistors on a chip doubles approximately every two years. While a larger number of transistors increases the processor's power, which is measured in floating point operations per second (FLOP/s), the increased transistor density leads to greater power dissipation.

Power dissipation is not only dependent on the density of transistors but is also a function of the third power of processor clock frequency [44]. Therefore, packing many processors with a smaller clock frequency onto a CPU significantly reduces the energy dissipation while being capable of holding the theoretical CPU's power constant.

As a result, the developmental frequency of new CPUs has been leveling flat in recent years, as is displayed in Figure 1. Instead, the increased number of logical cores keeps transistor counts going.

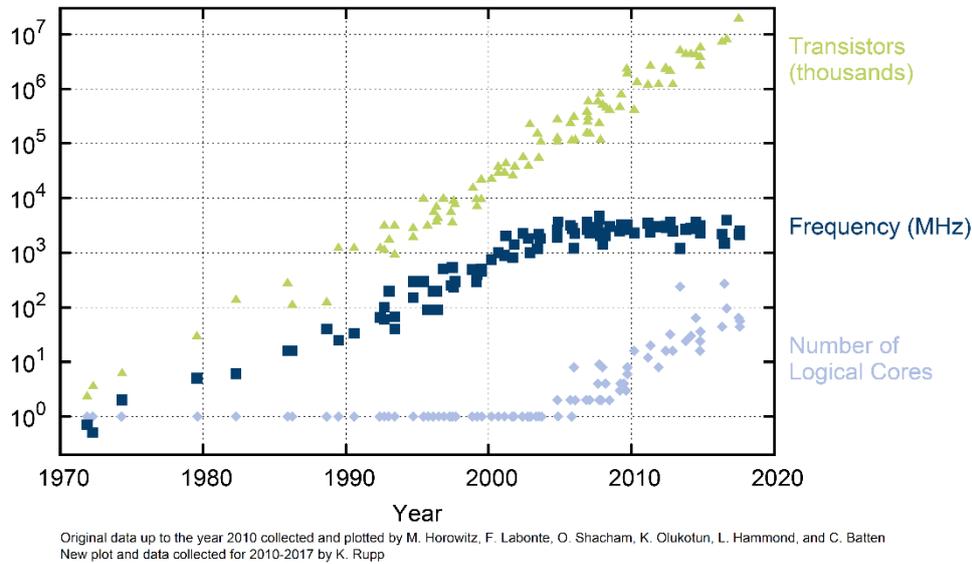

*Figure 1. Development of transistor counts, frequency, and number of logical cores, based on Rupp [45].*

**1.4 Computational limitations of optimization solvers**

To take advantage of theoretical CPU power, this in return leads to the necessity of decomposing the optimization problem into discrete parts that can be either processed independently or which must communicate minimal data between themselves.

Nevertheless, the majority of the commercially-available optimization solvers are not capable of exploiting these parallelized resources on larger scales. While an increase from one to four CPUs enables a significant reduction in the computational runtime [46], Rehfeldt et al. [47] show that in some instances of large-scale energy system models, the reduction of computational runtime stagnates at around eight computational threads. Similar observations have also been made for the general benchmarking of Mixed Integer Linear Programs, where no computational runtime improvement can be observed from 4 to 12 threads [48].

Therefore, today's supercomputers with a huge amount of cores, e.g., JUWELS with 122,768 cores or Sunway TaihuLight with 10,649,600, are not capable of efficiently tackling this mathematical complexity, and the advancements in computational resources cannot be accessed for ESOMs.

**1.5 Objective and structure of the paper**

As a result of the constraint exploitation potential of computational resources, methods to efficiently tackle and reduce the complexity of ESOMs are being continuously developed by the energy system research community. These can be either simple qualitative evaluations to identify irrelevant parts of the system models or sophisticated machine learning methods that systematically reduce the data input of energy system models.

Historically, many different complexity management methods have been individually introduced and benchmarked, with first attempts to compare different complexity reduction methods with respect to accuracy and computational impact, either for dispatch [49] or transmission expansion models [50].

Nevertheless, to the knowledge of the authors, a holistic and comprehensive overview of complexity management methods does not exist, and therefore this work reviews, evaluates, and qualitatively compares the methods to each other. Thereby, we want to guide modelers who encounter computational limitations and identify research gaps to lay a foundation for the future development of methods to reduce complexity.

The article is structured according to the process of system modeling: First, Section 2 provides an overview of what constitutes energy system optimization models and describes their complexity. Then, Section 3 discusses methods to systematically reduce or manage the complexity of energy system models. Section 4 introduces possibilities to decompose the resulting mathematical optimization models. The main conclusions will then be drawn in Section 5.

## 2 Determinants of complexity

Complexity is not an end in itself, although we can sometimes perceive it differently in the research community. Thus, it is important to first understand what drives and determines the complexity of an energy system and its respective model.

Therefore, in Section 2.1, we define energy systems from a system-theoretical perspective, identify energy system modeling as the process of depicting these in Section 2.2, compare different approaches to defining and measuring their complexity in Section 2.3, and finally derive conclusions about the dimensions and drivers of complexity in energy system modeling in Section 2.4.

### 2.1 Definition of energy systems from a system-theoretical point of view

Bertalanffy [51] defines systems as complexes of interacting elements that cannot be described by their elements alone. Today, this property is referred to as *emergence*. Laughlin [52] specifies that by looking at the individual elements of systems, the overall system's behavior is in many cases not recognizable [52]. Meadows et al. [53] expand on existing definitions of systems with the dimension of a system's *purpose*. In a technical environment, systems are defined as a "set of interrelated elements considered in a defined context as a whole and separated from their environment" (DIN IEC 60050-351, DIN German Institute for Standardization, 2014). This settles the existence of a defined *system boundary* and *holism*.

The existence of system boundaries results in a system hierarchy: systems can be part of other systems. Basic elements form a first system, a so-called sub-system, which in turn is part of a larger system, a super-system [53] [54]. Skyttner [54] introduces a fairly rigorous and thus well transferable definition of different levels of hierarchy: *parts* build *units* that form *components*, which can be summarized as *modules* that in turn build a *subsystem* as part of a *system* that is embedded in a *macrosystem*. Figure 3 summarizes the findings of this definition and provides an overview of the hierarchical structure of different system levels.

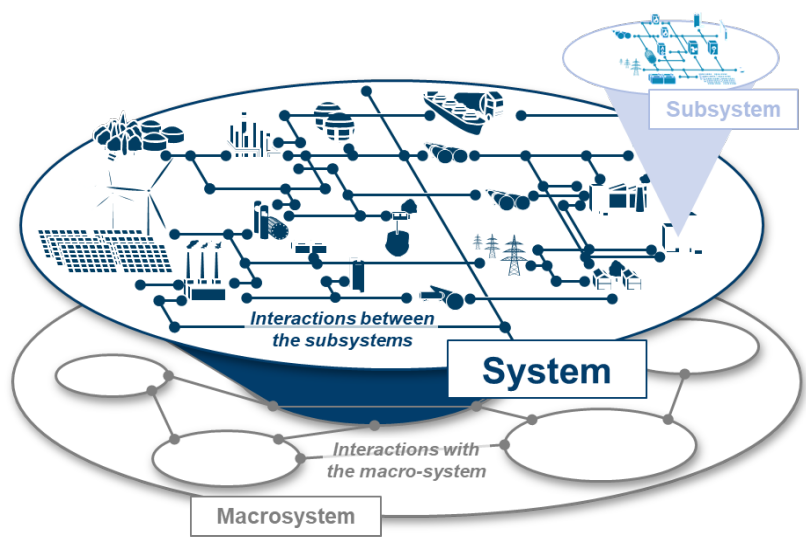

*Figure 2. Definition of an energy system in a hierarchy of system levels with their interaction. An example can be a building energy supply system as **system**, which is embedded into a district supply system as **macrosystem**, while single components such as a heatpump are **subsystems** with different technology components in themselves.*

Based on, energy systems can be defined as special types of systems. Table 1 demonstrates the key aspects of two sample energy systems.

*Table 1. Key aspects of two sample energy systems from a system-theoretical perspective.*

| **Energy system** | Heat supply of a building (A) | Energy supply of a country (B) |
|---|---|---|
| **System boundaries** | Building's exterior walls | State borders |
| **Purpose of the system** | Meeting the heating needs of building users | Meeting all energy needs of the country's population |
| **System elements (excerpt)** | *Technical components*: heat supply, pipes, heat storage(s), heatpump, walls, windows, doors, etc. *Agents:* inhabitants, guests, etc. | *Technical components*: power plants, electrical grid, storage, etc. *Agents*: consumers, producers, financiers, etc. |

Using the terminology introduced, we can state that system B serves as a macrosystem to system A, and vice versa system A is a subsystem of system B. Additionally, we can see that both systems comprise subsystems, modules, and components.

**2.2 Models to depict energy systems**

As experiments on energy systems themselves are not possible, or only possible at a very high cost, energy system models are frequently used tools to analyze system behavior [55], while different types and categorizations of energy system models can be defined [6, 56, 57]. In the scientific context, energy system models must *be purposeful, repeatable, unbiased,* and make a *novel contribution* [58].

In general, a model is a simplified representation of reality for a specific purpose. It, therefore, possesses three central properties [59]:

1. *Representation:* a model replicates a part of reality.
2. *Simplification:* a model does not replicate all properties of the original. Rather, it is designed to replicate only those that the modeler deems relevant.
3. *Pragmatism:* a model has a purpose dictating what part of reality it represents and how this is achieved.

Overall, energy system models try to serve a purpose, i.e., derive a specific decision by representing the most relevant parts of the analyzed systems and introducing necessary simplifications. Thus, the model seeks to depict the emergent behavior of the energy system under investigation, i.e., output variables and their dependence on external circumstances in the form of input variables, as is illustrated in Figure 4.

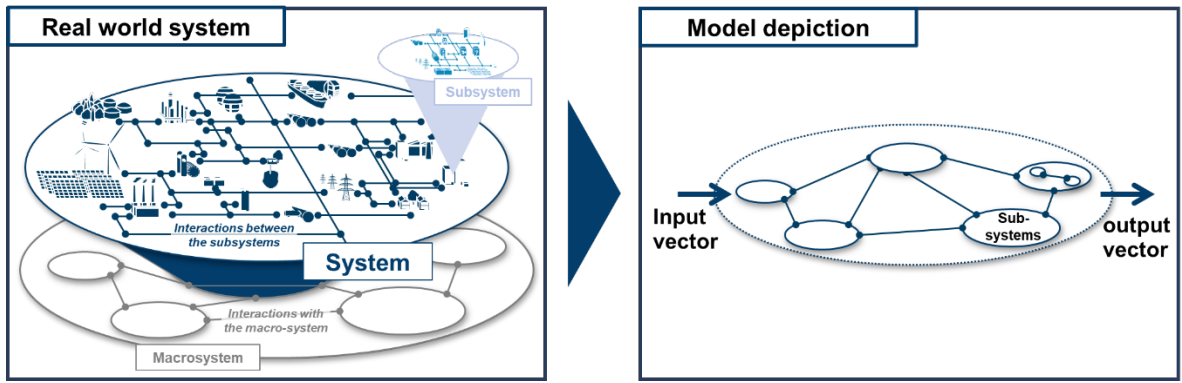

*Figure 3. Illustration of the modeling process as a selected and simplified representation of the reality. For the example of designing a building supply system, the possible grid supply options such as a heating network of the district system as macrosystem need to be considered in a model depiction at the system boundary, while certain subsystems such as the heatpump are simplified to single elements with a described behavior, e.g., the coefficient of performance.*

## 2.3 Complexity in energy systems and energy system models

In the following, we assess what drives the complexity of the original system and its simplifying model, whereby the definitions of measures for complexity are introduced.

In general, concepts for a formal description of complexity have been derived from system theory [60], mathematics, and the natural sciences [61-64]. More recently, there have been various attempts to define complexity [65, 66]. To this day, however, no concise and comprehensive definition of the concept has been established. Rather, definitions of complexity are highly context-dependent and driven by certain use cases. Given the sheer amount of different definitions and measures, Lloyd [67] opened a 'non-exhaustive list' that comprises over 40 different ways to describe complexity. This accords to the findings of Mitchell [68] (p. 301) that show that researchers focusing on complexity widely agree that we cannot yet characterize the phenomenon of complexity in a rigorous way.

According to the physicist Lloyd [67], three prevailing questions govern scientific concepts of complexity:

- How hard is it to describe the system?
- What is its inherent degree of organization?
- How hard is it to create a model?

The first two principles are much in line with the main characteristics of complexity, as described in the economic management literature [69]. Here, complexity is characterized by *variety* (i.e., the amount and kinds of elements in a system), *connectivity* (i.e., the amount and kinds of relationships between the elements), and *dynamics* (i.e., unpredictability). To handle complexity, three basic concepts are outlined in the management literature [70]:

1. *Complexity reduction*: reduce existing complexity by reducing the number of existing parts, variants, and processes, i.e., the number of variables and interdependencies in the system.
2. *Complexity control*: efficient control of unavoidable complexity through appropriate measures, i.e., appropriate anticipation of dynamics.
3. *Complexity avoidance*: avoid an increase in complexity beyond what is necessary.

It becomes clear that the definitions of complexity are very well transferrable to energy systems *per se*, whereas the measures to manage complexity mainly refer to energy system *models* and thus to the third question by Lloyd [67]. Therefore, we must distinguish the complexity of the underlying energy system and modeling complexity.

### 2.3.1 Complexity in energy systems

Based on complex system theory, the main characteristics of complex systems can be transferred to energy systems [58, 68, 71]:

1. The existence of agents in the system;
2. Networks that link different agents and physical components;
3. Dynamics in the sense of changes in time that might include feedback mechanisms;
4. Self-organization, meaning autonomous adaptation towards external changes;
5. Path-dependency, including lock-in effects;
6. Emergence that describes an emerging behavior in the system's macro structure;
7. Co-evolution regarding co-existence and interdependence with other systems;
8. Learning and adaption based on experimentation that leads to improved functionality of the system.

This list, in combination with the overview on definitions and measures for complexity, allows two main conclusions to be drawn: first, energy systems must be regarded as complex systems; second, the prevailing complexity of energy systems is currently increasing as their variety, connectivity, and dynamics grow.

### 2.3.2 Complexity in energy system models

According to Billings [72], the identification of relevant systems to model, i.e., the definition of system boundaries and relevant interdependencies between input and output variables, is a part of system theory. As is introduced above, models simplify and depict real systems. Thus, we conclude that as the complexity of energy systems increases, the complexity of energy system models that are intended to reflect their emergent behavior must also increase in turn.

Focusing on energy system models, the third question raised by Lloyd [67] becomes particularly important: how hard is it to create a model? The way this question is formulated leads to the conclusion that not only must the process of solving the model be considered but also the entire process, starting with collecting the necessary input data, followed by defining, implementing, and running the model and finishing with an interpretation of results.

Based on our synthesis regarding system theory and complexity research, we can now distinguish three different levels of complexity that are relevant to energy system modeling:

1. The complexity of the energy system itself (the left side of Figure 4);
2. The complexity of the part of the energy system to be investigated, i.e., the scope of the analysis (center of Figure 4); and
3. The resulting computational complexity of the model (right side of Figure 4).

While it is hard to find ubiquitous and quantitative measures for the first two types of complexity, it is more straightforward to measure the emergent computational complexity of the model. For the latter, we can distinguish two different approaches: hardware-dependent measures such as runtime and memory usage emergent in the level of complexity [73], or hardware-independent measures for algorithmic complexity in the sense of *Big-O*-notations, also referred to as *Bachmann-Landau* notation or asymptotic notation [74]. This measures the complexity of a given problem by means of the fastest formulation of a solution algorithm based on computer models such as the (deterministic) Turing Machine (Turing, 1937). Using this notation, algorithms can be classified by complexity [75, 76]. These classes provide upper limits to the increase in algorithmic complexity with the size of the input data, as is illustrated in Table 1.

*Table 2. Algorithmic complexity classes to describe how runtime or memory growth with the size of the input data (excerpt).*

|  | Constant | Logarithmic | Linear | Polynomial | Exponential |
| --- | --- | --- | --- | --- | --- |

| Big O notation | $O(1)$ | $O(\log n)$ | $O(n)$ | $O(n^a)$ | $O(2^n)$ |

Furthermore, algorithms can be divided into the categories of efficiently solvable (e.g., $O(n^2)$) and inefficiently solvable (e.g., $O(2^n)$) problem classes [77, 78], where the latter class is often referred to as 'intractable' or NP-hard problems. Another class is defined as 'advantageously parallelizable' problems [79]. The solution to problems of this class can be significantly simplified by implementing a shared memory. Overall, the algorithmic complexity definition has the advantage of being independent of constantly improving hardware performance. However, this only provides upper limits for the scalability of the algorithm with the size of the input data and therefore does not allow for valid comparisons across the complexity classes.

### 2.4 Accuracy and complexity in energy system modeling

As a result of the insights summarized above and in accordance with the findings by De Carolis et al. [80], we can conclude that there exists a trade-off between accuracy and complexity within energy system models: the more the modeler abstracts from the real system, the more the model strays from reality and the less its behavior reflects the real system's complexity. However, system theory also shows that it is sometimes sufficient to reflect the emergent behavior of complex systems and it is not always necessary, or sometimes even impossible, to reproduce this behavior by depicting all of a system's elements [52]. We can further conclude that modelers must indeed *resist the temptation to concentrate exclusively on computational complexity*. It is also important to bear in mind the underlying system's complexity drivers, as they will ultimately govern the complexity of the analysis. Here, the choice of system boundaries, depth of modeling, and the level of detail of the depicted interdependencies determine the complexity and accuracy of the model.

The trade-off between complexity and accuracy has been analyzed, e.g., by Bale et al. [81], leading to the conclusion that it also influences the communicability of results to a non-scientific audience (for instance policy-makers). Additionally, studies from different fields of research have recently focused on the trade-off between the accuracy and complexity of models [49] and did *not find a general superiority of more complex models*, as in cases of overfitting in the case of poor-quality input data [82, 83]. Another direction of research has focused on identifying optimal solutions within the trade-off between computation time (as a proxy for complexity) and model accuracy [84]. Finally, Brooks and Tobias [85] claim that there is a scarcity of research focusing on choosing the best-fitting models in terms of performance to depict the underlying system. They introduce four measures to quantify model performance: the quality of results, the future usability of the model, the verification and validation of model results, and the required resources.

It is thus imperative to identify those areas in which the trade-off is particularly favorable and those elements that considerably increase the complexity of a model without markedly increasing the accuracy of its results. In other words, it is necessary to: (1) identify the complexity drivers of energy system models; and (2) choose the right level of complexity for the present scope of analysis. For the latter task, a tension arises between the claim for complexity by Stirling [86] and the call for reductionism and '*parsimony*' in energy system models by De Carolis et al. [80].

## 3 Methods for complexity reduction

While the complexity of the energy system itself cannot be influenced, and that of the modeling process depends on the expertise of the modeler, *computational complexity* can be systematically altered and quantified.

The dimensions to reduce complexity of ESOMs are illustrated in Section 3.1, while the subsequent sections (Section 3.2, Section 3.3, Section 3.4, and Section 3.5) describe possibilities for reducing these.

### 3.1 Dimensions to reduce complexity

As discussed in Section 2, a general definition of complexity dimensions is challenging. From a mathematical point of view, however, three basic factors have a direct impact on computational complexity, as is illustrated in Figure 5: The total size of the optimization model, the optimization problem class, and the connectivity within the model.

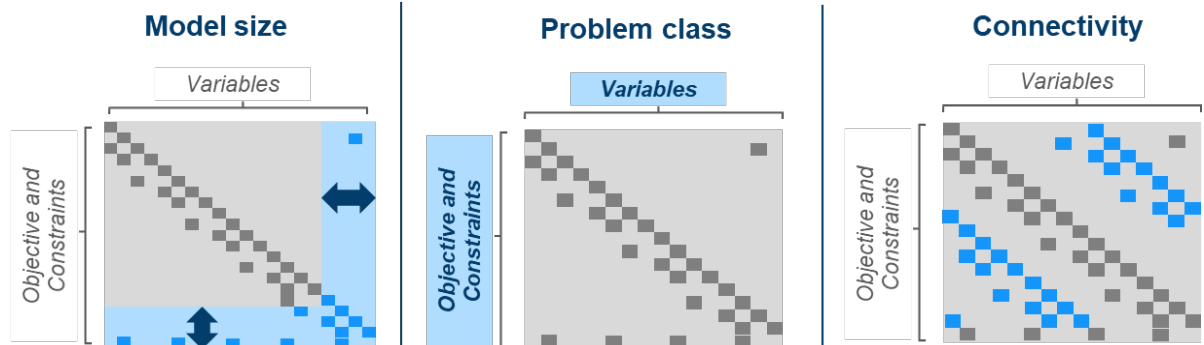

*Figure 4. Complexity dimensions expressed by sketching the objective and constraint matrix. The 'model size' depicts the general number of variables and constraints. The 'problem class' is derived from the type of variables (continuous, binary, or integer) and the type of constraints (linear, nonlinear). The 'connectivity' describes the linkage between the variables and constraints.*

The scope, spatial resolution, and temporal resolution of the ESOM determine the **model size** and thus the size of the matrix as a combination of constraints and variables: It can be influenced by the number of considered time steps, sometimes referred to as time slices or snapshots [38], in operation or investment decisions to account for their dynamic behavior. Its cardinality is determined by the temporal resolution, e.g., sub-minutely [87], sub-hourly [88, 89], or, in most cases, hourly [18], and the observation period can range from a number of typical days to a time series spanning decades [37]. We discuss the possibilities for reducing the temporal scale in Section 3.2. Likewise, the spatial and technological resolution directly affects the model size and, respectively, the size of the network. In particular, a high number of spatially-distributed nodes is required for optimal grid design [30, 37, 90, 91], but also for an adequate consideration of renewable supply technologies [92]. The options to lower the spatial scale of ESOMs will be framed in Section 3.3. Nevertheless, even single node ESOMs can result in large-scale optimization problems if a high resolution of technology types and sectors is considered [93, 94], as if technical solutions for heating, cooling, electricity, and industrial processes are regarded simultaneously, resulting in large networks. However, the latter strongly depend on the system's scope and the research question, wherefore no systematic complexity reduction methods currently exist to the authors' knowledge.

To achieve large-scale ESOMs, Linear Programs (LPs) are primarily used as a **problem class** [37, 40], which combine linear constraints with a continuous variable set. Their convex nature combines with historically significant efforts to develop efficient solving algorithms with polynomial solving time. As discussed above, they are especially required to design large-scale, bottom-up energy systems supplied by renewable power [14]. Recent works have also shown that convex Quadratic Programs (QP) are suitable for this problem scale [95]. Specific Mixed-Integer Linear Programs (MILPs) with a small amount of binary or integer variables for the technology choices [30, 90, 91, 96-99] allow for the optimization of larger networks for a few time steps. Yet, every binary or integer variable leads to a cut in the solution space, resulting in a non-convexity and an np-hard problem. Nonlinear performance functions, e.g., the part-load efficiency of a fuel cell, determine a non-convex set of operational states, and with it, a non-convex optimization problem, usually resulting in Mixed-Integer Nonlinear Programs (MINLPs). This problem class is computationally-intensive, which limits the size of the considered systems and/or the temporal observation time. Often, they are simplified to a MILP by either modelers [100] or by the solving algorithm itself [101]. The systematic simplification of technology models by avoiding nonlinearities or discontinuities and the related non-convexity of the program is habitually performed by modelers and will be discussed in Section 3.4.

The **connectivity** of an energy system model is indicated by the density of the constraint matrix for the linear case. An energy system with strong spatio-temporal linking expressed by dense transmission networks, operational dynamics such as states of storages, or investment dynamics, leads to high connectivity. For example, the modeling of a hierarchical order of time grids [102] [90] [103] [104] leads to a small model size but has a strong linkage within the model, which makes it hard to decouple parts of it, making it challenging to solve different parts within a parallelized computer infrastructure. For single-core performance, this is challenging to generalize: As Tejada-Arango et al. [105] show, implementing operational characteristics comes with a trade-off between compactness (i.e., increasing the model size) and tightness (i.e., the similarity between the solution search space of the MILP and the relaxed LP) [106]. They find that ramping constraints in particular positively impact the solving time by increasing the tightness. Otherwise, as in the case of dispatch models, the weak connectivity is often neglected by temporally decoupling the model in favor of being able to solve a large-scale network for different time steps in a parallel computer infrastructure [83]. For a more limited model scope, however, it is possible to include dynamic time-coupling [107]. The options to tackle this connectivity from a modeler's perspective are discussed in Section 3.5 and lead to the possibility of solving and decomposing ESOMs, as described in Section 4.

With the help of this definition of the computational complexity of ESOMs, we can now describe how the original complexity drivers of energy systems, published by Bale et al. [58], relate to the computation of their virtual representatives in Table 1.

*Table 1. A holistic list of complexity drivers in energy systems and their assumed impact in computational representations*

| Complexity drivers in energy systems | Drivers for computational complexity |
| --- | --- |
| **Agents** in the system | Increased **model size** (more nodes to depict) and more complex **problem class** to depict behavior |
| **Networks** that link different agents and physical components | Increased **model size** (more connections to depict) |
| **Dynamics** in the sense of changes in time that might include feedback mechanisms | Increased model **connectivity** and increased **model size** to resolve those dynamics, maybe even the need for a more complex **problem class** due to non-linearities |
| **Self-organization**, meaning autonomous adaption to external changes | Increased **connectivity** of the model, perhaps even the need for a more complex **problem class** due to non-linearities, **model size** due to complex behavior |
| **Path dependency** including lock-in effects | Increased **connectivity** of the model and the need for a more complex **problem class** due to non-convexity |
| **Emergence** behavior of the system's macrostructure | More complex **problem class** due to non-linearities |
| **Co-evolution** regarding co-existence and interdependence with other systems | Increased **model size** (need to depict interactions with macrosystems) |
| **Learning and adaption** based on experimentation that leads to improved system functionality | Increased model **dynamics** and a more complex **problem class** due to non-linearities |

### 3.2 Temporal aggregation

As introduced above, the temporal complexity of ESOMs is spanned by the temporal resolution and the considered time horizon [38]:

- A long time horizon is crucial for the generation of expansion-planning models to appropriately capture transformation pathways [6], technical learning rates [108], or other long-term effects such as economic, social, or environmental processes that gain importance at different time scales [109].
- In contrast, a high temporal resolution is crucial for considering the rising share of highly intermittent renewable energy sources, e.g., in cost-optimal unit commitments [18, 110].

The basic idea of temporal aggregation is to represent the time series of demands, supplies, or residual loads with a large number of time steps by a smaller number of time steps [38]. This is achieved either by the direct reduction of the temporal resolution, as discussed in Section 3.2.1, or by representing time series with a reduced number of typical periods, explained in Section 3.2.2. The general concept of these two approaches is highlighted in Figure 6.

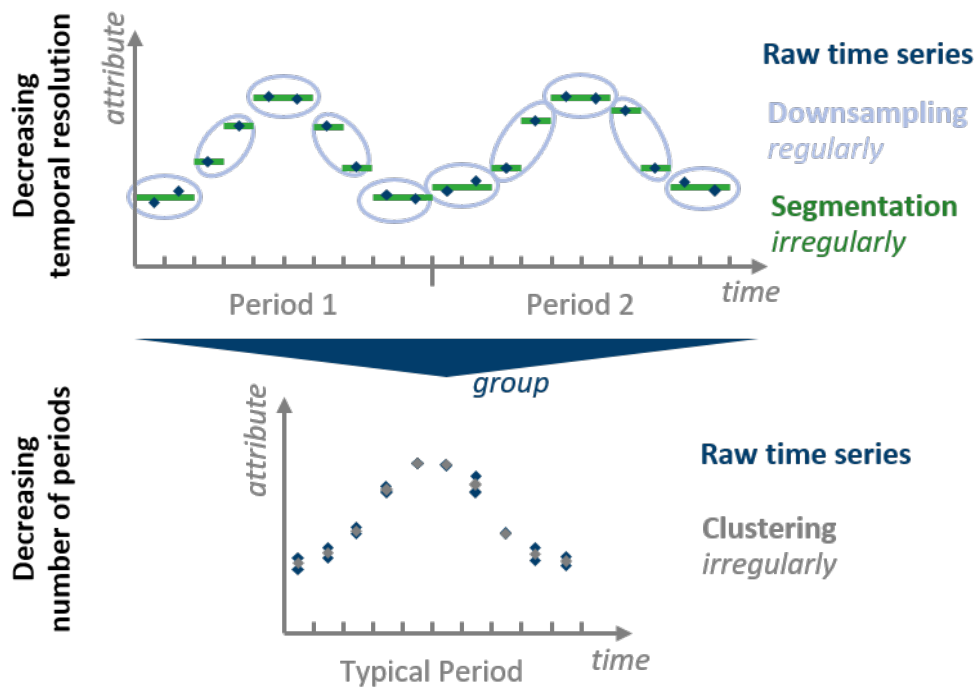

*Figure 5. Overview of the possibilities for reducing the temporal complexity in an discretized time space*

### 3.2.1 Decreasing the temporal resolution

The simplest way to reduce the overall number of time steps is to decrease the temporal resolution directly. As is illustrated in the upper part of Figure 5, this can either be done in a regular manner, referred to as down-sampling, or in an irregular manner, based on the similarity of adjacent time steps, called segmentation.

For down-sampling, a predefined number of adjacent time steps is taken and represented by its mean, e.g., by taking the average of every two time steps. This leads to an underestimation of the time series' variance of the raw input data and can have a severe impact on the result of energy system optimization models, such as for systems with a high share of renewable energy. Here the necessary capacities to be built are underestimated [37], or the self-consumption rate of energy systems, including feed-in from photovoltaic and battery storage are overestimated [36, 87], and also counts for feed-in from wind turbines [110].

A more advanced method is segmentation, the merging of adjacent time steps based on their mutual similarity. This approach leads to irregular new time step lengths, as highlighted by the green bars in the upper part of Figure 5. Numerous approaches exist [38]: Some favor partitional clustering [111] or agglomerative clustering [112] under the constraint that only adjacent time steps are to be clustered, while other approaches include

MILP optimizations to merge adjacent time steps while minimizing the deviation from the original time series [113].

### 3.2.2 Decreasing the number of periods

The concept of typical periods is that not only might adjacent time steps within a series be similar, but also entire periods within the time series. This is illustrated in the lower part of Figure 5, where periods 1 and 2 are merged based on their comparable profiles.

The similarity between different periods is achieved by cutting normalized input time series into typical periods and aligning them in a dissimilarity matrix, which is illustrated in Figure 7. Each row can be interpreted as a hyper-dimensional candidate point for the clustering procedure.

Clustering algorithms generally strive to maximize the intra-cluster similarity of candidate points while maximizing inter-cluster dissimilarity by allocating data points to different clusters [114]. In the realm of time series clustering, this means that "homogeneous time series data are grouped based on a certain similarity measure" [115]. The cluster centers or selected data points are then chosen as representatives or "typical periods" [116-118].

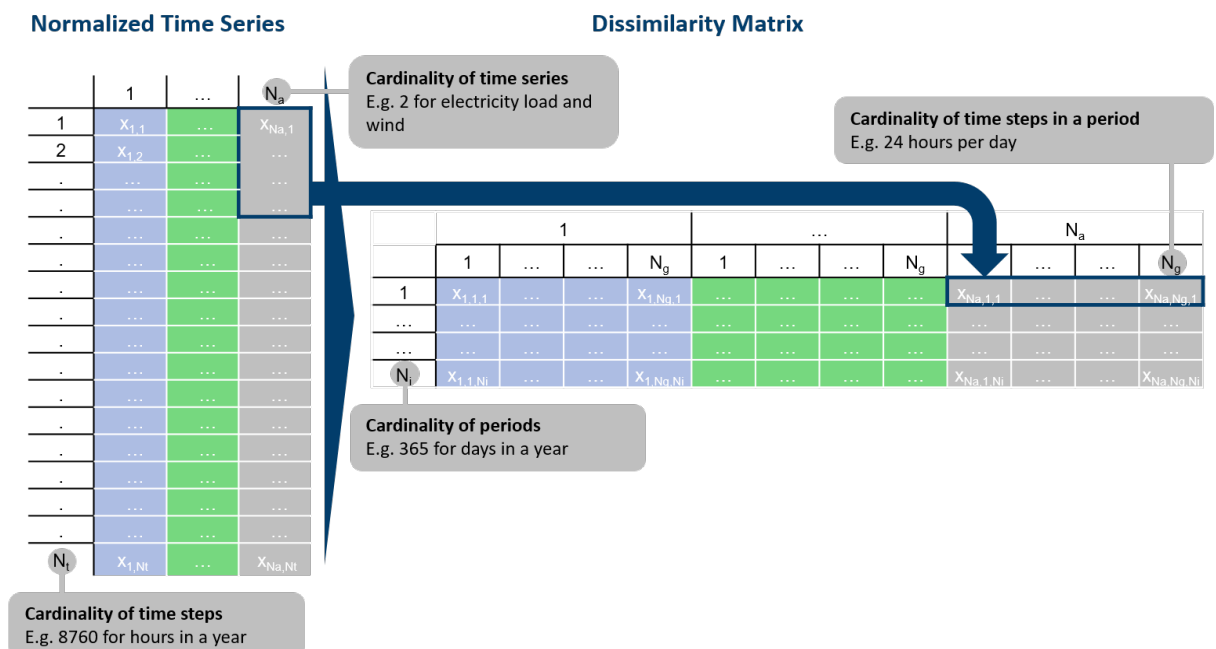

*Figure 6. Dimensionality reduction step to transform a set of normalized time series to candidate periods in a dissimilarity matrix which is used as basis for the clustering process.*

While one of the first cluster method applied was fuzzy c-means clustering [119], the most common cluster methods use Euclidian distance as a metric, k-means or Ward's [120], hierarchical greedy clustering algorithms with either medoids (e.g., existing typical periods) or means (e.g., synthesized typical periods) as representatives, or the MILP formulation of the k-medoid algorithm [121], which searches for representative days that minimize the distance to all of the other data points of the respective cluster [116, 118, 122, 123].

Recently, new clustering algorithms were introduced and tested on energy system models that take temporal shifts into account for a certain bias along the time axis, such as the k-shape algorithm, which was introduced [124] and applied to the thermal energy demands of university buildings [125] and electricity prices [126].

Other publications have focused on daily duration curves as candidates for clustering and removing intra-day order [30, 127]. Another method is to directly approximate certain values in a yearly duration curve and find a combination of daily duration curves whose linear combination minimizes the error to the original duration curve [35, 128, 129]. This method is also capable of taking the fluctuation and variance of duration curves into account [130]. However, these methods widely neglect the intra-daily order of typical days and therefore may fail to model the intra-day dynamics.

Moreover, the issue of an appropriate number of clusters is the subject of current research. Although a vast number of clustering indicators [115] such as, for example, the silhouette score [131, 132] and distortion sum ratios between homogeneously distributed data samples and real samples [133, 134] exist and have been applied as indicators for a sufficient number of typical periods [111, 135], no set of indicators has proven to be superior over all of the others. A reason for this could be that many input time series describe continuous phenomena, which means that the sample points are not well-separated from each other, and so the clustering error simply monotonically decreases with an increasing number of typical periods [111, 136].

Last but not least, it is to highlight that clustered periods do not respect the chronology of typical periods and therefore require alternative description of dynamics between the periods, e.g., for the consideration of seasonal storages. Some solutions exist [104, 137] but they are challenging to implement and can alter the results of the analysis [138] wherefore proper validations are required.

### 3.2.3 Extreme periods

The design of energy systems strongly depends on extreme periods, which are crucial to operational feasibility and surplus capacities, as well as the operational costs for which the prevalence of the different typical periods is important [139].

A common method is to define days containing the peak or minimum value of a certain attribute as an extreme period, e.g., maximum demand days or minimum supply days with respect to solar or wind energy production [117, 118, 140]. Moreover, other authors have noted that not only are single extreme values of interest for the robustness of optimizations, but also cumulative extreme values, such as days with cumulative peak demand [141, 142]. Apart from that, typical days can also be considered an extreme period if they contain the peak or minimal demand within a season [143]. Finally, it is worth noting that other heuristics have also been proposed that focus on maintaining a certain variance or gradient within a clustered time series [122] or by directly importing them into the optimization problem of the clustering algorithm [139].

Recently, other methods have proposed tackling the problem of underestimating the variance of time series with respect to clustering in general. Some add iteratively feasible time steps if the operation of an energy system optimized for an aggregated time series is not feasible for operation with the original time series [144] while others use synthesized variations [18, 142, 145] in time series and simulate the operation of energy systems designed for one scenario with all of the others [142], or re-run the optimization, including the most expensive time steps from the first optimization run [145, 146].

Moreover, aggregation methods have been introduced that define the upper and lower bounds of the original optimization problem and try to arrive at the solution by iteratively increasing the number of time steps [147-150].

All things considered, three main directions can be identified in the development of more robust results from energy system optimizations: Heuristic approaches that focus on preserving the important characteristics of the original time series by directly preserving it in the aggregated time series, searching for the most robust energy system by creating multiple scenarios and systematic approaches based on aggregated input data to define the upper and lower bounds for the objective of fully-resolved energy system optimization

### *3.3* **Spatial aggregation**

The spatial complexity of energy system models is determined by two aspects, namely the spatial resolution, hence the number of defined model regions, and the representation of these, i.e., the number of technologies or agents per model region. Fundamentally, we define model regions as geometries containing energy system components. These regions are connected to other regions by inter-regional connections, referred to as a network or grid. Additionally, the spatially-distributed energy system components within regions are connected. However, these intra-regional connections are generally assumed to be copper plates, such that the data of all energy components for each region are generally simply aggregated [50]. Given different research questions and the resulting requirements concerning the level-of-detail of these energy systems, spatial resolution and representation varies significantly and thus requires the spatial aggregation of the available energy system data.

Spatial aggregation generally comprises both the *grouping* of regions with similar properties on the one side [151] and the *representation* of the regions' information, such as time series within defined regions, using aggregate functions on the other [152], as is seen in Figure 8. Thus, the grouping determines how to aggregate the network, whereas the representation determines how to aggregate the technologies within the newly created regions. Therefore, aggregation functions are required to represent the data of the initial set of regions for the newly reduced region set. These aggregate functions range from simple representations with a sum, mean, min, or max [153], to more advanced functions such as median, mode, and rank [153], to complex aggregation functions in the form of algorithms, such as network reduction algorithms, as described by Hörsch and Brown [154].

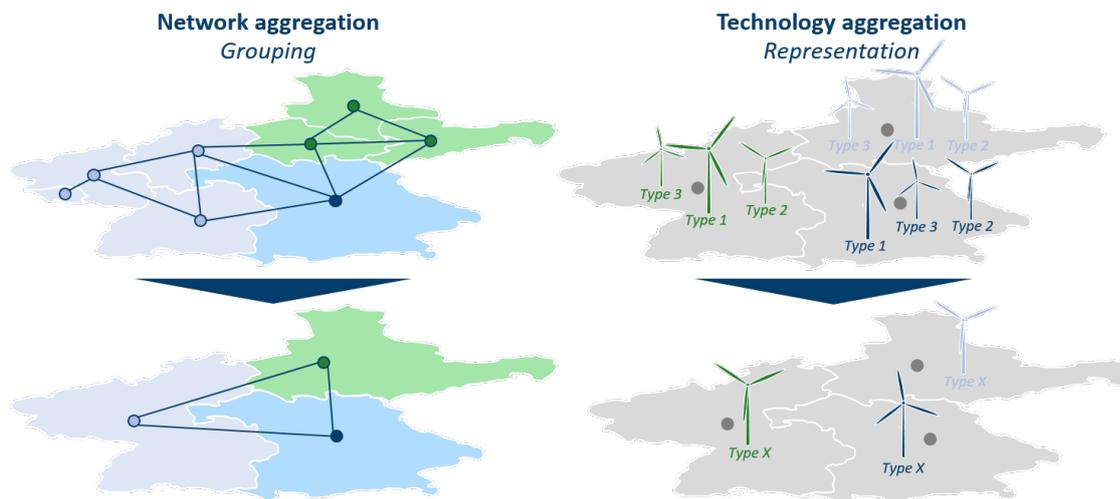

*Figure 7. Difference between the aggregation of a network of nodes and different technologies inside a modeled node.*

So far, the grouping was primarily based on administrative boundaries, with current energy systems analysis studies with a European scope mostly being based on national regions, with national time series of hourly temporal resolution [32]. The scope of these studies ranges from electricity-only scenarios [155] to sector-coupled ones [92, 156, 157].

Initial studies have been conducted to analyze the impact of spatial aggregation on energy systems analysis based on grid constraints [154, 158, 159]. Whereas Hörsch and Brown [154] derive a clustered equivalent network with the help of k-means based on the geometric distance of load and conventional generation, Svendsen [159] clusters based on calculated power transfer distribution factors, while Fazlollahi et al. [160] cluster based on demand patterns.

Furthermore, Scaramuzzino et al. [161] analyzed a grouping based on NUTS3 regions using multiple indicators to identify similar regions in terms of energy potential indicators, such as onshore wind energy and photovoltaics, as well as non-energy related indicators, such as GDP and population density. On the other hand, Siala et al. [162] analyzed groupings based on highly resolved renewable energy potential and demand data. On the de-

*Submitted to Advances in Applied Energy* 14

mand side, clustering algorithms to represent similar spatially-distributed buildings [38] or industrial sites can lead to higher quality energy demand representations. Another example is the aggregation of municipalities [163] that is used to determine a representative set of municipalities that can each obtain its individual energy supply design. Figure 9 summarizes the different existing spatial aggregation applications in different layers but also frames the research gap of defining spatial aggregation in terms of a holistic consideration of all layers that are relevant for energy system design and operation.

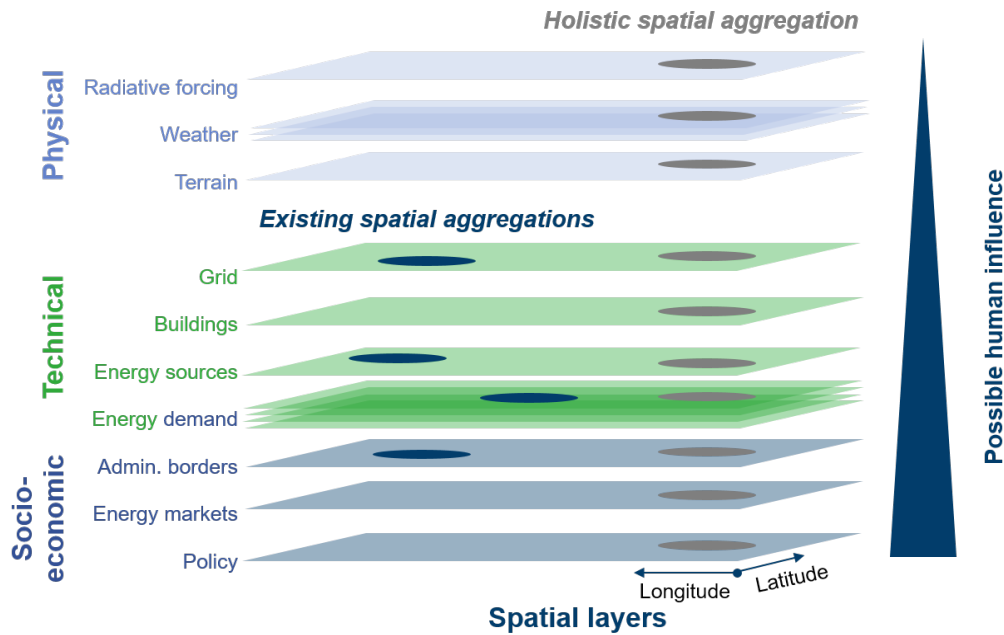

Figure 8. Illustration of different spatial layers that should be considered in the spatial definition of energy system models.

The impact of spatial aggregation on computational complexity was noted by Cao et al. [50]. The results show a significant dependence on the number of regions for both computational indicators and the obtained solution in terms of aggregated costs. This motivates further investigation of spatial aggregation techniques to determine an appropriate number and selection of regions as a function of the complexity of the subsequent energy system optimization. Nevertheless, it is important to distinguish between an aggregation in an integrated ESOMs context [154, 162], in contrast to an aggregation of independent entities that become separately optimized [163, 164]. For the latter, only a linear computational complexity reduction can be achieved relative to the aggregation rate, but their representation can be more easily achieved, as the connection between the candidates does not need to be considered.

In contrast, for connected systems, aggregation results in information loss due to representation inside every region, *and* because of the balancing effects between different regions that can lead to an underestimation of data variabilities. In particular, the copper plate assumption for intra-regional connections results in an externalization of costs. These could be either explicitly additionally accounted for [160] or neglected [154]. An alternative would be to internalize them using respective representation functions. On the supply side, the spatial distribution of renewable energy sources with many different turbine designs – on- or offshore, optimized for full-load hours or maximal generation - and Photovoltaic system configurations – orientation and inclination – require sufficient representation as well, for example by clustering similar capacity factor time series inside every grouped region. All in all, the information loss that is relevant for the energy system optimization must be systematically identified and minimized.

## 3.4 Reduction of the level of detail in modeling system behavior

The mathematical description of the system operation determines the resulting problem class of the ESOM. While accurate modeling would in general lead to a MINLP, the technical characteristics – that in reality would limit flexibility in the operation or the choice of the technical components – are simplified or omitted, leading to an MILP or even an LP.

Some of the general simplifications can be generalized for ESOMs and are grouped into four categories:

1) Constraints that are *continuous nonlinear* (e.g., specific investment costs that decrease with the unit size) can be linearized using different linearization methods that are not mathematically equivalent to the original equation but rather an approximation with a linear formulation [101, 165].

2) *Non-continuous* constraints, such as either-or, if-then-else constraints, absolute value or minimum value functions can be linearly implemented in optimization models using the so-called Big-M method [166-169] and natively result in an MILP, as the Big-M method adds a binary variable and a sufficiently large arbitrary value *M* to the model.

3) *Multi-dimensional functions* of the form $f(x_1,x_2)$ (e.g., efficiencies that depend on the realized size of a generation unit and the operational load level) can be linearly implemented by fixing all but one variable to predefined values and precomputing the results [165].

4) Finally, products that *multiply a continuous with an integer* variable can be substituted by a continuous variable and post-processed to the original variables [170].

Table 2 displays a summary of some linearization methods and the resulting problem formulation for the four types of nonlinear constraints, while their application is discussed in the following section for operation and investment modeling in ESOMs.

*Table 2. Categories of nonlinear and non-continuous constraints in energy system optimization models and respective linearization methods*

| Type | Example in ESMs | Linearization Method | Resulting Problem | Assessment of the computational burden | References |
|---|---|---|---|---|---|
| Continuous | Specific investment costs | Binary Steps / SOS Type 1 | MILP | High: Adds multiple binary variables | [165] |
| | | Piecewise-Linear Function / SOS Type 2 | MILP | High: Adds multiple binary variables | [165] [101] |
| | | Intercept slope | MILP | Medium: Adds a single binary variable | [171] |
| | | Constant Value | LP | Low: Nonlinearities are omitted and no binary variables are added | [172] |
| Non-continuous | Minimum loads | Big-M | MILP | Medium: Adds a single binary variable | [173] |

| Nested function | Efficiency depends on the load and unit size | Precomputing with binary steps | MILP | Low: Simplified representation by the discretization of continuous non-linear relationships without adding binary variables | [165] |
| --- | --- | --- | --- | --- | --- |
| Integer product | Total load of unit group | Substitution | MILP | Low: Simplified representation using discretization without adding binary variables | [170] |

We discuss the implications of operation modeling in Section 3.4.1 and investment modeling in Section 3.4.2.

### 3.4.1 Level of detail in the operation of system components

Energy conversion, transmission, or storage units underlie several technical restrictions. Common technical restrictions applied in energy system models include part-load-dependent efficiencies (PLDE), minimum part loads, start-up and shut-down costs, minimum down- and up-times, as well as ramping rates [49]. The PLDE is the relationship between the output and input energy during an energy conversion process, as illustrated in Figure 10. The PLDE depends on the load level of the conversion unit, making the relationship a continuous nonlinear function.

This nonlinearity can be approximated using various methods. The use of a single constant value (e.g., the PLDE at full load) is based on the simplified assumption that no part load dependency exists, resulting in a simple LP formulation. This is the most trivial but also least complex way of linearizing PLDE, and is commonly performed in energy systems analysis [174, 175]. Binary steps and piecewise-linear approximations introduce binary variables [176]. An intercept-slope only adds one binary variable per unit to the model; however, it can only be applied to certain bounded functions. The more sections the original function is divided into the more accurate and the more complex the approximation becomes. Special ordered sets (SOS) of Type 1 or Type 2 can be used to implement binary step or piecewise linear approximations [177]. New approaches consider even piecewise quadratic approximations of the PLDE, resulting in Mixed Integer Quadratic Constrained Programs (MIQCP), but could not be proven to outperform conventional MILP formulations [178].

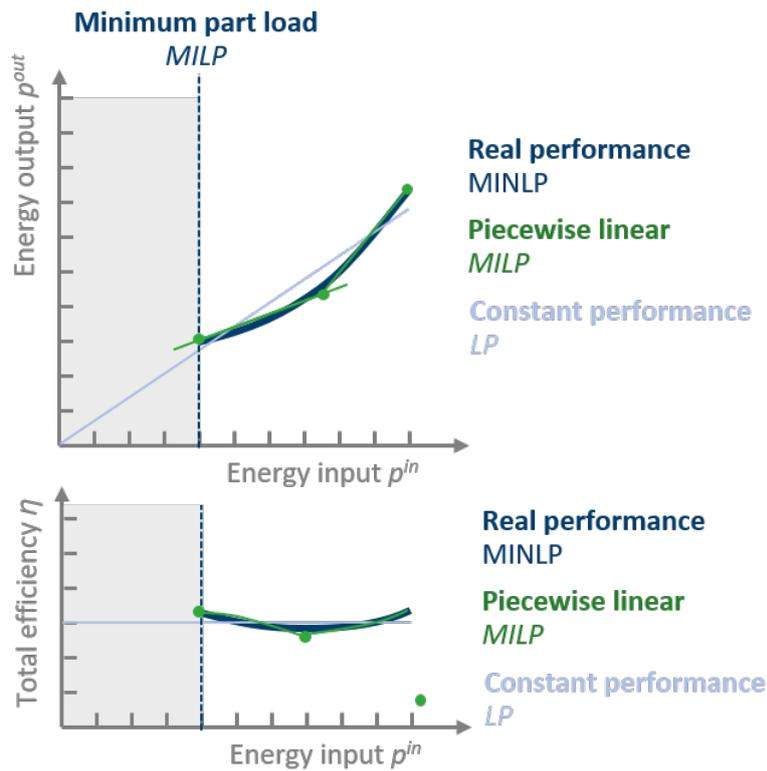

*Figure 9. Possibilities for modeling the input-output performance operation of an energy system technology and the resulting mathematical program type.*

Conversion units usually only run above a minimum load, e.g., for turbines or engines of 20% to 50% of their rated power [179]. This distinction of the load being above or below the minimum load is modeled with a binary variable [173] in combination with the Big-M method and is also illustrated in Figure 10.

The starting up or shutting down of a unit can cause additional fuel and depreciation costs [179]. To add these costs to the cost function of an optimization problem, binary variables for the unit status *start-up* and *shut-down* must be introduced in addition to the aforementioned binary variables for the on and off statuses. To set the cost, the information on the change in the *operating state* between two consecutive time steps is required [180]. This is realized by connecting the binary variables through dynamic constraints. Some units must remain shut down for a minimum downtime or kept running for a minimum uptime. This additional required information on the unit's status from previous time-steps entails a dynamic constraint, as is shown in [106] and illustrated in Figure 11.

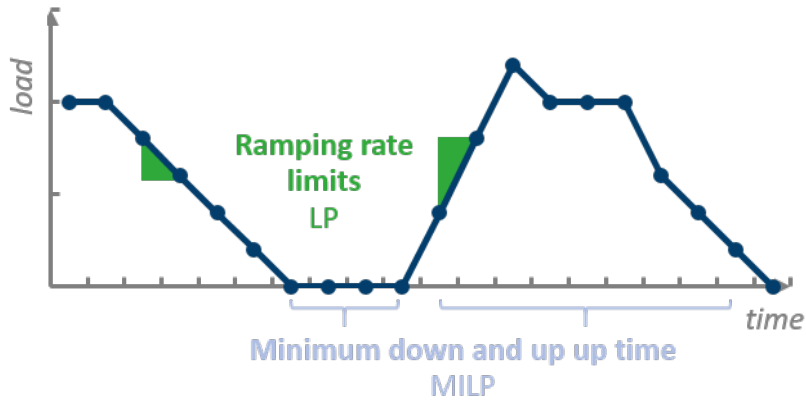

*Figure 10. Typical dynamic constraints of a technology operation and their mathematical program type.*

The rate with which a unit alters its load can be restricted by upper and lower ramping limits due to technological constraints. In addition, ramping up and down can impose costs for which an additional continuous variable is introduced that is added to the cost function [181]. Restricting the load alternation between consecutive time steps requires dynamic constraints. Costs for ramping up or down require the already introduced binary variables in the case where they are to be distinguished from start-up and shut-down costs.

In addition to the technical complexity, with an increasing amount of deregulated electricity markets, there is also an increased effort in mathematical modeling of related sociatal systems, see, e.g., [182], such as competitive equilibrium in electricity markets to predict the market clearing price (MCP) and respecting its related uncertainty in systems operation [183].

**3.4.2    Level of detail in investment cost modeling**

Additional to the operational constraints of the systems, Figure 12 illustrates the challenge to model the choice and scaling of the technologies.

The choice of technology unit available on the market, including their related price and performance, would introduce many single *binary or integer* variables [100, 144, 184, 185]. This is computationally challenging, but closest to reality. Especially in small-scale energy system models, in which only a small discrete number of components has to be chosen for achieving cost-optimality, discrete investment decisions, which is also referred to as "lumpy investments" [20, 186] can have a considerable impact on the cost-optimality. Further, this also applies to investments with fixed and size-independent cost contributions, e.g. development expenses.

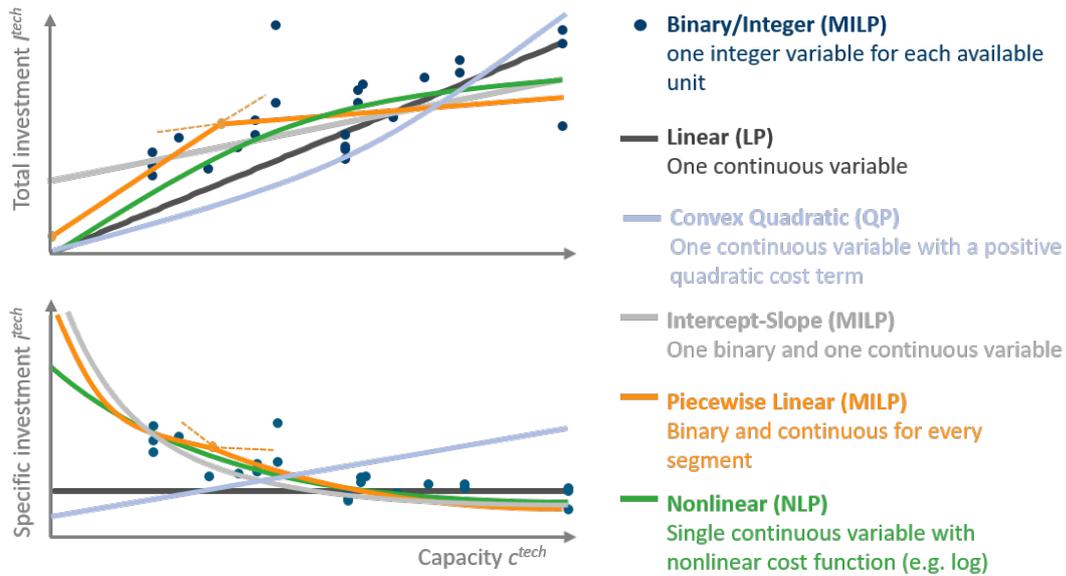

*Figure 11. Constituted mathematical program type depending on the chosen investment cost model.*

Aggregated energy system models on a global or national scale do not require this detail, as they rely on abstract perspectives of the system and often model technology scaling with continuous *linear* cost functions [40, 93, 187, 188].

Nevertheless, the consideration of learning effects in macroeconomic models determines *nonlinear* cost curves that reduce the specific costs of a technology with an increase in its capacity. These are especially common in national or global energy assessment models with small temporal resolutions [189]. To combine those with bottom-up approaches, this non-linearity can be linearly approximated piecewise [170]. In contrast, another possibility for aggregated energy system models is to introduce a *convex quadratic* cost term that increases the specific investment cost for higher capacities [95]. The idea is that small capacities can be exploited at low cost, e.g., wind turbine locations, but the higher the capacities become, the more challenging and cost-intensive will be the deployment. An advantage is that it does not significantly increase the computational complexity in comparison to a linear approach due to its convex solutions space.

For the design of local energy systems, it is common to approximate the technology investment with a cost share related to their existence (*intercept*) and cost share related to the scale (*slope*) [25, 99, 190, 191], resulting in an MILP. The choice of efficiency measures, e.g., in the building envelope, constitutes binary variables [25, 192]. While the Intercept-Slope approach already respects economies of scale for a small range, it still has high estimation errors for larger ranges of the technology scale, such as CHP units. Those *nonlinear* cost functions require also *piecewise linear* approximations [141, 165, 167, 193], generating more binary variables but providing a sufficient degree of accuracy.

### 3.5 Simplification of system dynamics and connectivity

Another option to manipulate the complexity of the model is by neglecting or reducing its dynamics or intrinsic connectivity between different decision variables.

One approach to decreasing the connectivity of the model components is to separate the decisions for the technology choice, sizing and operation, as is visualized in Figure 13. This can either be done fully separated in an iterative manner, such that in the top layer a genetic optimization algorithm defines the design of the ener-

*Submitted to Advances in Applied Energy*      20

gy system, and a bottom layer wherein either a simulation [93, 192, 194, 195] or separate operation optimization [196-200] defines the operation.

Another alternative is to simplify the operation model in the design layers, e.g., by fixing the impedance in transmission expansion models in the first stage [201, 202] or by aggregating the time resolution of the operational model to determine a design and then validating and fine-tuning it with full temporal resolution [127, 149, 203].

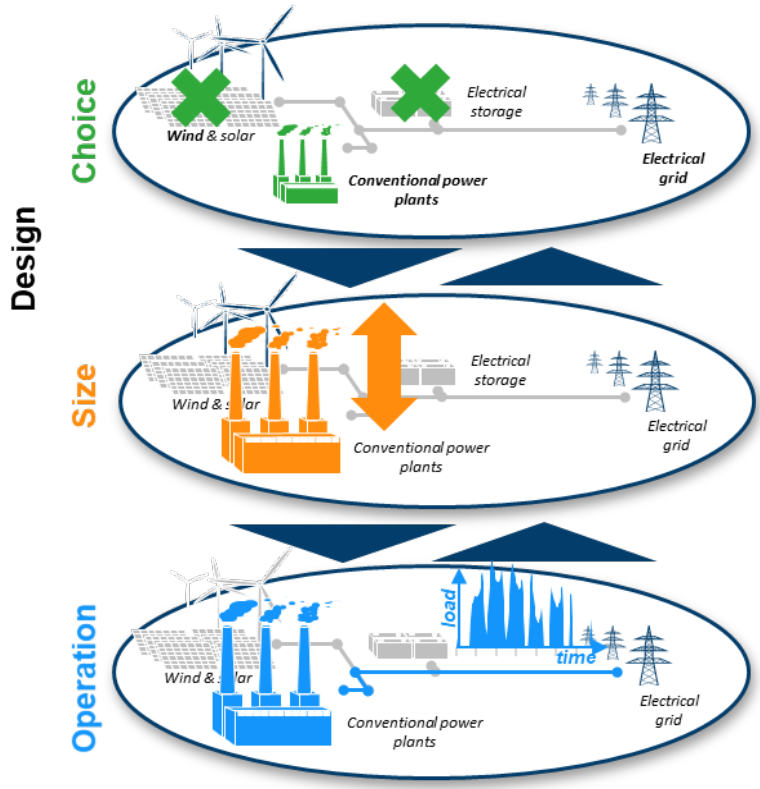

*Figure 12. Decision layers influencing the total energy system design.*

Another possibility is to separate decisions within each of the layers. For the case of capacity expansion, i.e., models applied for determining the design of the energy system, the connectivity between consecutive transformation phases can be implemented using two fundamentally different modeling approaches: A perfect-foresight and a myopic approach, the latter also being known to be limited or restricted foresight [6]. These are illustrated in Figure 14, but can also be transferred to operation modeling.

Perfect-foresight is based on complete information about the past and future requirements for the energy system model. This means that the past and future constraints of all expansion phases in the mathematical model are known at any time step. Accordingly, perfect-foresight approaches are capable of finding a cost-minimal transformation pathway across all expansion phases [189]. In contrast, myopic approaches assume limited knowledge about the future, meaning that optimization in each expansion phase is based on the results of the previous expansion phases and the constraints of the current expansion phase only. In that way, myopic models more appropriately capture short-sighted decision-making under real economic conditions [6, 189, 204].

| Connectivity | Run # | 2020 | 2030 | 2040 | 2050 |
|---|---|---|---|---|---|
| **No foresight** | 1 | ● | | | |
|  | 2 | | ● | | |
|  | 3 | | | ● | |
|  | 4 | | | | ● |
| **Myopic** Foresight | 1 | ● | | | |
|  | 2 | | ● | | |
|  | 3 | | | ● | |
|  | 4 | | | | ● |
| **Myopic** Back-casting | 1 | | | | ● |
|  | 2 | ● | | | |
|  | 3 | | | ● | |
|  | 4 | | | ● | |
| **Rolling horizon** | 1 | ● | ● | | |
|  | 2 | | ● | ● | |
|  | 3 | | | ● | ● |
| **Perfect Foresight** | 1 | ● | ● | ● | ● |

*Figure 13. Classification of modeling approaches for capturing expansion and investment dynamics across a long time horizon based on Lichtenböhmer et al. [205] and Lopion et al. [6].*

Apart from the extreme cases of approaches, i.e., strictly perfect foresight models and completely myopic models, a number of mixed approaches also exist for energy system design. One example is the rolling horizon [188, 206-208]. A rolling horizon splits the full observed interval into several smaller intervals and matches the first value of an interval with the last value of the previous iteration. Thereby, the entire time horizon is split into expansion phases, but in such a way that the temporal subsets overlap. Thus, the overlapping time subsets are linked to each other in a myopic manner, while the optimization within each subset is run with perfect foresight, which enables the model to consider limited foresight [205] while drastically reducing the computational burden compared to a perfect-foresight approach [206]. This raises the question, however, of optimal interval lengths for the trade-off between complexity and accuracy.

The same applies to operation modeling, where models are distinguished into models with decoupled time steps, rolling horizon models, and perfect foresight models. Models with decoupled time steps implement small operational characteristics of technical components, as many of such characteristics would result in coupled time steps. These models are commonly used to assess network loads [143]. Models applying rolling horizons can implement most of the operational characteristics of technical components mentioned above. The main difference to models using perfect foresight is the representation of long-term storage. While approaches that use a rolling horizon must pass storage values from one period to the next (cf. [208]), approaches rely on a perfect foresight model for long-term storage with a consistent variable.

Nevertheless, the separation of the different model elements has recently consisted of manual decomposition purposely performed by modelers in order to reduce the complexity or natively due to a limited model scope. Although these methods will converge on a solution for their underlying ESOMs, it is not natively global and must be subsequently validated. In consequence, exact decomposition methods that also quantify the error could be advantageous.

## 4 Solving and Decomposition methods

It is important to know what makes a problem difficult to solve, despite the use of modern computing hardware.

Until the 1990s, optimization problems were divided into linear and non-linear types. This trend changed due to significant work by Rockafeller [209] and optimization problems are therefore now viewed as convex or non-convex. Non-convex problems are significantly harder to solve than convex ones; this is because any local optimal solution to a convex problem is guaranteed to also be a globally optimal solution. When sufficient computing resources are available, modelers can use decomposition methods to exploit parallel computing for both convex and non-convex optimization classes.

**4.1 Solving and decomposition methods for convex optimization**

Typical convex optimization problems are linear or certain types of quadratic problems; the standard algorithms for these two classes are the simplex algorithm or interior point methods and the simple interior point methods, respectively.

Linear programs can ordinarily be very efficiently solved by commercial solvers. Nevertheless, for large scale problems it can become necessary to break the original problem into smaller parts that are coupled via a so-called coupling or the linking of variables and constraints. This is called decomposition while the most common decomposition methods for convex optimization problems are Lagrangian relaxation [210], Benders decomposition [211], and Alternating Direction Method of Multipliers (ADMM) [212]. When implemented in parallel [213, 214] significant reductions in the compute time can result [215].

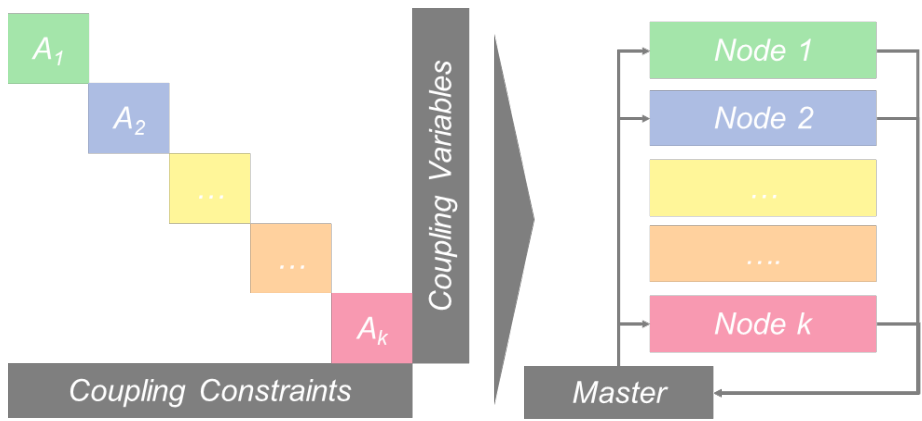

*Figure 14: Sketch of the arrow head shape of the constraint matrix for an efficienct decomposition.*

Although the idea exist, to automatically decompose optimization programs, an efficient decomposition requires that the work be divided into, at best, equally-sized work packages. Finding an optimal decomposition is generally not possible in polynomial time [216, 217], and so in solving these problems one resorts to heuristics and approximations. An example is the Generic Column Generation [218] that automatically detects those structures by comprising an arrowhead and bordered detector via graph partitioning using *hmetis* [219], as shown in Figure 15. Other graph partitioning frameworks are those of *Chaco* [220] and *Scotch* [221]. This decomposed problem may then be fed into a structure-exploiting algorithm such as Parallel Interior Point Solver (PIPS) [47, 222]. However, generating this block structure automatically takes for large problems significant compute time, wherefore a the structure identification is most often done by the modelers themselves.

One example is the temporal decomposition, where for each time segment we obtain a smaller problem, as is shown in Section 3.5. The linking variables and constraints ensure that the states at the end of one time segment are the same as those at the beginning of the next.

As a further example for the usage of decomposition methods within energy system optimization, Stursberg [223] applied the Benders decomposition to the convex formulated capacity expansion problem with an improved approach for cut generation.

## 4.2 Solving and decomposition methods for non-convex optimization

The presence of integer variables in an optimization problem is typically enough to make it non-convex. Non-convex problems can be either MIPs or MINLPs. MIPs are solved by iteratively solving multiple LPs, while MINLPs are first converted into either multiple NLPs or MIPs that are then broken down further, for instance as proposed by Goderbauer et al. [100] for the case of an energy system. More generally, we iteratively relax (by "branching") and enforce constraints (by "bounding") on the original problem, the first algorithms of this type dating back to the 1960s [224-226]. Relaxations enlarge the feasible set of the MIP without excluding the feasible points of the original problem. In the best case scenario, they provide a problem that is easier to solve and for which an optimal solution can be obtained faster than the true problem, which allows the user to derive a lower (upper) bound for a minimization (maximization) objective. The efficiency of branching procedures depends on two strategic decisions concerning (i) the selection of the branching variable; and (ii) the selection of the next node that must be solved. Several studies attend to these decisions and evaluate them (e.g., [227-229]), while this theoretical work is well-reflected in commercial solvers.

Decomposition methods can be employed within branch and bound algorithms as well. The generalized Benders decomposition method proposed by Benders [211] and extended by Geoffrion [230] with an outer approximation algorithm by Duran & Grossmann [231] has been an immensely popular method to date for solving optimization problems. As an additional step in the branch-and-bound algorithms, it is also possible to generate so-called cutting planes and add them as new constraints in order to cut off infeasible solutions [101, 232].

In the context of energy system optimization, decomposition methods such as nested Benders decomposition are mainly used in stochastic scenario analyses intended to account for the uncertainties in power generation unit commitment [233], wind power investment decisions [234], or sustainable energy hub designs [235]. Other use cases of decomposition methods are the solving of large-scaled, multi-period MILP problems in the context of electric power infrastructure planning and optimal power flow taking transmission constraints into account, see, e.g., [236-238]. More examples of the decomposition methods for energy optimization can be found in [239].

## 5 Conclusions

Based on the review of complexity drivers in energy systems optimization models (ESOM), we derive a qualitative guide to support modelers in their modeling.

Complexity goes beyond computational complexity: The entire process of including the necessary input data, followed by defining, implementing, and running the model and interpreting the results is complex and must be communicated. The plain result "42" is not sufficient. Therefore, **consciously define the ESOM's superstructure** and start with a coarse, simple model to identify the necessary determinants to answer your research question, and then increase the level of detail where necessary. In general, a model design based on data availability should be avoided and instead choosen the research goal.

A high temporal-spatial resolution and temporal-spatial scope gains importance with a higher share of renewable energy sources and directly impacts the size of the ESOM and its calculation time. Therefore, a **systematic reduction of the size of the model** is recommended. Simple down-sampling approaches on the temporal scale allow one to quickly identify the relevance of the chosen temporal resolution. If this is given, superior time series aggregation by clustering typical periods are gaining popularity but can require significant adaptations to the model formulation, e.g., for respecting the chronology of periods. Similar counts for the spatial resolution: Aggregating based on neighboring administrative regions can result in quick computational gains; nevertheless, the impact of the underlying copper plate assumptions inside each region must be evaluated. Holistic aggregation schemes simultaneously based on grid structures, renewable potential, and demand are challenging to implement due to the heterogeneous data types but will be relevant in future to evaluate the spatial scope.

Technical and economic relationships are natively non-linear or non-continuous. Still, the majority of the reviewed ESOMs are linear and continuous, **resulting in convex optimization models** that have polynomial solving times. In consequence, binary variables should be avoided and equeatiosn linearized where possible. For instance, a technology with a negligible cost contribution does not need a piecewise linear cost function. Mixed-integer linear models or non-convex non-linear programs are np-hard, constituting an exponential solving time which strongly limits their application for renewable energies in combination with storage and transmission technologies.

In order to still achieve feasible system designs and operations, the **errors must be quantified**. This can either be done by a benchmark model or by more advanced error-bounding and multi-stage approaches that evaluate the error due to aggregation or model simplifications. As conservative estimators, upper bounds can even guarantee feasibility of the original problem. Nevertheless, the latter are challenging to implement and often specific to the simplification.

Advancements in computational resources have lately been achieved through an increase of the number of cores while the calculation frequency stagnates. In consequence, modern software must be parallelized, meaning that the optimization model must be decomposed. This has until now not been automatically possible with any free or commercial optimization solver at larger scales. In consequence, exact **decomposition methods must be tailored to the model** which is time consuming while the computational improvements are strongly related to the connectivity of your model. Heuristic decompositions can do the job as well but are not necessarily globally optimal.

Last, based on this guide and review, some major research gaps were identified:

- A holistically-quantified cross-impact analysis of complexity reduction methods, including model simplifications, aggregation methods, and heuristic decompositions for different energy models is open. In particular, the impact of systematic avoidance binaries and non-linearities in larger energy system models must be quantified.
- To the authors' knowledge, no methods that aggregate the entire optimization model based on its abstract mathematical description exist, although such methods could significantly accelerate the finding of start solutions in optimization solvers or support the bounding process.
- Lastly, improvements in the solving algorithms are required in order to more easily exploit modern parallel computing infrastructure.

# 6 Acknowledgments

The authors acknowledge the financial support of the Federal Ministry for Economic Affairs and Energy of Germany for the project METIS (project number 03ET4064).

# 7 Authors Contribution

Conceptualization, L.K., M.R.; Validation, A.P., D.S., M.R.; Investigation, L.K, L.N., M.H., T.G., A.S., J.P., H.B., R.B., F.K., B.S.; Writing –Original Draft, L.K, L.N., M.H., T.G., A.S., J.P., H.B., R.B., F.K., B.S.; Writing –Review & Editing, L.K, M.H., B.S., M.R.; Funding Acquisition, A.P., D.S., M.R.; Resources, A.P., D.S.; Supervision, L.K., A.P., D.S., M.R.

# 8 Appendix

*Table 3: Overview to some bottom-up energy systems optimization frameworks. For a more detailed overview of their capabilities, it is referred to Groissböck [31], Lopion et al. [6], and Ringkjøb et al. [32].*

| Name | Type | Applications | Reference |
|---|---|---|---|
| LEAP - Long-range Energy Alternatives Planning | Annual time-step simulation for a mid- to long-term system planning horizon. | Different scales: Cities, states, national (e.g., Phil- | [15] |

| | | ippines) and global applications. | |
|---|---|---|---|
| EFOM - Energy Flow Optimization Model | LP of a national or multinational technology network to supply the consumer demand. | European Union with each country as single node. | [16] |
| BESOM - Brookhaven energy system optimization model | Single time step dispatch and planning optimization of a technology network to supply electric and non-electric demands. | National energy system of the United States. | [12] |
| MARKAL - MARKet and Allocation | Optimization of dispatch and long term planning of energy systems with a set of representative time slices. | Different scales: Community, state, national and global. | [17] |
| MESSAGE - Model for Energy Supply Strategy Alternatives and their General Environmental Impact | Gams-based LP used for capacity expansion planning and scenario analyses | A historically grown model used for different applications from capacity expansion to macroeconomic and climate impacts of energy supply and demand, | [14, 18] |
| IKARUS - Instrumente für Kilmagas-Reduktionsstrategien | LP optimization for planning and dispatch of a national supply system with with myopic-foresight and an operation in typical days. | National supply system of Germany. | [6] |
| PERSEUS - Programme-package for Emission Reduction Strategies in Energy Use and Supply-Certificate Trading | GAMS based dispatch optimization. | Germany and Europe. | [19] |
| TIMES - The Integrated MARKAL-EFOM System | GAMS based LP modeling framework based on time slices for energy system design. | International and national energy systems around the globe with adaptions to local systems. | [20-23] |
| DESOD - Distributed Energy System Optimal Design | C# based MILP used for district planning. | Capacity expansion planning of an energy system for a residential and commercial district comprising heat and electricity. | [24] |
| DER-CAM - Distributed Energy Resources Customer Adoption Model | GAMS based MILP microgrid planning framework with a dispatch based on typical days. | Different microgrid systems around the globe. | [25] |
| CALLIOPE (No acronym, c.f. [32]) | Python based LP modeling framework. | Applied to Europe, Great Britain, Italy, South Africa, China, Kenya, Cambridge and Bangalore. | [26] |
| OEMOF - Open Energy Modeling Framework | Python based MILP modeling framework for planning and dispatch of energy systems. | Different regional energy systems from single sites to municipalities or states, especially located in Germany. | [27] |
| URBS (No acronym, c.f. [32]) | Python based LP framework for capacity expansion planning and unit commitment. | Its usage has been adapted from neighbourhoods to continents. | [28] |
| PYPSA - Python for Power System Analysis | Python based software for LP optimization and power flow simulation of large scale power system networks. | Power grids across the globe while the its application to the European power grid is probably the most | [29] |

| | | popular one. | |
|---|---|---|---|
| FINE - Framework for Integrated Energy System Assessment | Python based MILP energy systems modeling framework that can be applied to different spatial and temporal scales. | Applied to European and German heat, gas and electricity networks as well as district systems. | [30] |